# CONDITIONAL HAAR MEASURES ON CLASSICAL COMPACT GROUPS

By P. Bourgade

*Université Paris 6*

We give a probabilistic proof of the Weyl integration formula on $U(n)$, the unitary group with dimension $n$. This relies on a suitable definition of Haar measures conditioned to the existence of a stable subspace with any given dimension $p$. The developed method leads to the following result: for this conditional measure, writing $Z_U^{(p)}$ for the first nonzero derivative of the characteristic polynomial at 1,

$$\frac{Z_U^{(p)}}{p!} \stackrel{\text{law}}{=} \prod_{\ell=1}^{n-p}(1-X_\ell),$$

the $X_\ell$'s being explicit independent random variables. This implies a central limit theorem for $\log Z_U^{(p)}$ and asymptotics for the density of $Z_U^{(p)}$ near 0. Similar limit theorems are given for the orthogonal and symplectic groups, relying on results of Killip and Nenciu.

**1. Introduction, statements of results.** Let $\mu_{U(n)}$ be the Haar measure on $U(n)$, the unitary group over $\mathbb{C}^n$. By the Weyl integration formula, for any continuous class function $f$,

$$(1.1) \quad \mathbb{E}_{\mu_{U(n)}}(f(u)) = \frac{1}{n!}\int\cdots\int |\Delta(e^{i\theta_1},\ldots,e^{i\theta_n})|^2 f(e^{i\theta_1},\ldots,e^{i\theta_n})\frac{d\theta_1}{2\pi}\cdots\frac{d\theta_n}{2\pi},$$

where $\Delta$ denotes the Vandermonde determinant. Classical proofs of this density of the eigenvalues make use of the theory of Lie groups (see, e.g., [6]), raising the question of a more probabilistic proof of it.

PROBLEM 1. Can we give a probabilistic proof of (1.1)?

This problem will be shown to be closely related to a suitable definition of conditional Haar measures. From the Weyl integration formula, the measure of $(e^{i\theta_1},\ldots,e^{i\theta_{n-p}})$ conditionally on $e^{i\theta_{n-p+1}} = \cdots = e^{i\theta_n} = 1$ $(1 \leq p \leq n)$ is (we omit the normalization constant in the following)

$$(1.2) \qquad \prod_{1 \leq k < l \leq n-p} |e^{i\theta_k} - e^{i\theta_l}|^2 \prod_{j=1}^{n-p} |1 - e^{i\theta_j}|^{2p} d\theta_1 \cdots d\theta_{n-p}.$$

We are interested in the converse problem of finding a measure on $U(n)$ with the eigenvalues statistics (1.2).

PROBLEM 2. What is the Haar measure on $U(n)$ conditionally on the existence of a stable subspace of dimension $p$?

We will give a natural matrix model, based on a product of independent reflections, for this conditional expectation, inducing the measure (1.2) on the spectrum.

Under this Haar measure on $U(n)$ conditioned to the existence of a stable subspace of fixed dimension $p < n$, let $Z_U^{(p)}$ be the first (i.e., the $p$th) nonzero derivative of the characteristic polynomial at 1. In the specific case $p = 0$ (i.e., under the Haar measure), the distribution of $Z_U^{(p)}$ has great importance: it allowed Keating and Snaith [8] to conjecture the moments of the Riemann zeta function along the critical axis, supporting Hilbert and Pólya's idea that there may be a spectral interpretation of the zeroes of $\zeta$. In a recent work [3], it was shown that

$$(1.3) \qquad Z_U^{(0)} \stackrel{\text{law}}{=} \prod_{k=1}^{n}(1 - X_k),$$

the $X_k$'s being independent with an explicit distribution. This result was generalized to other compact groups in [4].

PROBLEM 3. Under the Haar measure on $U(n)$ conditionally on the existence of a stable subspace of dimension $p$ is there an analogue of (1.3) for $Z_U^{(p)}$?

This decomposition exists and is a natural byproduct of the matrix model we give to answer Problem 2. It will allow us to give a central limit theorem for $\log Z_U^{(p)}$ and asymptotics of its density. We will also discuss the counterparts of these results for the symplectic and orthogonal groups: such asymptotics were related to the averages over families of elliptic curves in [13] and [14]. Our work is organized as follows.

*Conditional Haar measures.* Let $(e_1, \ldots, e_n)$ be an orthonormal basis of $\mathbb{C}^n$, and $r^{(k)}$ $(1 \leq k \leq n)$ independent reflections in $U(n)$ (i.e., $\mathrm{Id} - r^{(k)}$ has rank 0 or 1). Suppose that $r^{(k)}(e_j) = e_j$ for $1 \leq j \leq k-1$ and that $r^{(k)}(e_k)$ is uniformly distributed on the complex sphere $\{(0, \ldots, 0, x_k, \ldots, x_n) \mid |x_k|^2 + \cdots + |x_n|^2 = 1\}$. Results from [3] and [4] state that

$$r^{(1)} \cdots r^{(n)}$$

is Haar-distributed on $U(n)$. We generalize this in Section 2, showing that

$$(1.4) \qquad r^{(1)} \cdots r^{(n-p)}$$

properly defines the Haar measure on $U(n)$ conditioned to the existence of a stable subspace of dimension $p$.

*A probabilistic proof of the Weyl integration formula.* The previous conditioning of Haar measures allows to derive the Weyl integration formula (1.1) by induction on the dimension, as explained in Section 3. Our proof relies on an identity in law derived in the Appendix. As a corollary, we directly show that the random matrix (1.4) has the expected spectral law (1.2).

*The derivatives as products of independent random variables.* For the conditional Haar measure from Section 1 [$\mu_{U(n)}$ conditioned to the existence of a stable subspace of dimension $p$], we show that the first nonzero derivative of the characteristic polynomial is equal in law to a product of $n - p$ independent random variables:

$$(1.5) \qquad \frac{Z_U^{(p)}}{p!} \stackrel{\mathrm{law}}{=} \prod_{\ell=1}^{n-p} (1 - X_\ell).$$

This is a direct consequence of the decomposition (1.4). The analogous result for orthogonal and symplectic groups will be stated, directly relying on a previous work by Killip and Nenciu [9] about orthogonal polynomials on the unit circle (OPUC) and the law of the associated Verblunsky coefficients under Haar measure. Links between the theory of OPUC and reflections will be developed in a forthcoming paper [5].

*Limit theorems for the derivatives.* The central limit theorem proved by Keating and Snaith [8], and first conjectured by Costin and Lebowitz [7], gives the asymptotic distribution of $\log Z_U^{(0)}$. More precisely, the complex logarithm of $\det(\mathrm{Id} - xu)$ being continuously defined along $x \in (0, 1)$, they showed that

$$(1.6) \qquad \frac{\log Z_U^{(0)}}{\sqrt{1/2 \log n}} \xrightarrow{\mathrm{law}} \mathcal{N}_1 + i\mathcal{N}_2$$

as $n \to \infty$ with $\mathcal{N}_1$ and $\mathcal{N}_2$ independent standard normal variables. The analogue of their central limit theorem for the Riemann zeta function was previously known by Selberg [12] and can be stated as

$$\frac{1}{T} \int_0^T dt \mathbb{1}\left(\frac{\log \zeta(1/2+it)}{\sqrt{1/2 \log \log T}} \in \Gamma\right) \xrightarrow{T \to \infty} \int \int_\Gamma \frac{dx\, dy}{2\pi} e^{-(x^2+y^2)/2}$$

for any regular Borel set $\Gamma$ in $\mathbb{C}$. Section 4 explains how the identity in law (1.5) allows to generalize (1.6):

$$\frac{\log Z_U^{(p)} - p \log n}{\sqrt{1/2 \log n}} \xrightarrow{\text{law}} \mathcal{N}_1 + i\mathcal{N}_2.$$

From results in [9], a similar central limit theorem holds for the orthogonal and symplectic groups, jointly for the first nonzero derivative at 1 and $-1$. This will be related to works by Snaith ([13] and [14]) in number theory; her interest is in the asymptotic density of $Z_{SO}^{(2p)}$ near 0. Using identities in law such as (1.5), we give such density asymptotics for the classical compact groups.

**2. Conditional Haar measure.** For $r$ a $n \times n$ complex matrix, the subscript $r_{ij}$ stands for $\langle e_i, r(e_j) \rangle$, where $\langle x, y \rangle = \sum_{k=1}^n \overline{x}_k y_k$.

2.1. *Reflections.* Many distinct definitions of reflections on the unitary group exist, the most well known may be the Householder reflections. The transformations we need in this work are the following.

DEFINITION 2.1. An element $r$ in $U(n)$ will be referred to as a reflection if $r - \mathrm{Id}$ has rank 0 or 1.

The reflections can also be described in the following way. Let $\mathcal{M}(n)$ be the set of $n \times n$ complex matrices $m$ that can be written

$$(2.1) \qquad m = \left(m_1, e_2 - k\frac{\overline{m}_{12}}{1 - \overline{m}_{11}}, \ldots, e_n - k\frac{\overline{m}_{1n}}{1 - \overline{m}_{11}}\right)$$

with the vector $m_1 = (m_{11}, \ldots, m_{1,n})^\top \neq e_1$ on the $n$-dimensional unit complex sphere and $k = m_1 - e_1$. Then the reflections are exactly the elements

$$r = \begin{pmatrix} \mathrm{Id}_{k-1} & 0 \\ 0 & m \end{pmatrix}$$

with $m \in \mathcal{M}(n-k+1)$ for some $1 \leq k \leq n$. For fixed $k$, the set of these elements is noted $\mathcal{R}^{(k)}$. If the first column of $m$, $m_1$, is uniformly distributed on the unit complex sphere of dimension $n - k + 1$, it induces a measure on $\mathcal{R}^{(k)}$, noted $\nu^{(k)}$.

The nontrivial eigenvalue $e^{i\theta}$ of a reflection $r \in \mathcal{R}^{(k)}$ is

$$(2.2) \qquad e^{i\theta} = -\frac{1 - r_{kk}}{1 - \overline{r_{kk}}}.$$

A short proof of it comes from $e^{i\theta} = \text{Tr}(r) - (n-1)$. We see from (2.2) that for $r \sim \nu^{(k)}$ this eigenvalue is not uniformly distributed on the unit circle, and converges in law to $-1$ as $n \to \infty$.

2.2. *Haar measure as the law of a product of independent reflections.* The following two results are the starting point of this work: Theorems 2.2 and 2.3 bellow will allow us to properly define the Haar measure on $U(n)$ conditioned to have eigenvalues equal to 1. These results appear in [3] and [4], where their proofs can be found.

In the following we make use of this notation: if $u_1 \sim \mu^{(1)}$ and $u_2 \sim \mu^{(2)}$ are elements in $U(n)$, then $\mu_1 \times \mu_2$ stands for the law of $u_1 u_2$.

THEOREM 2.2. *Let $\mu_{U(n)}$ be the Haar measure on $U(n)$. Then*

$$\mu_{U(n)} = \nu^{(1)} \times \cdots \times \nu^{(n)}.$$

THEOREM 2.3. *Take $r^{(k)} \in \mathcal{R}^{(k)}$ ($1 \leq k \leq n$). Then*

$$\det(\text{Id} - r^{(1)} \cdots r^{(n)}) = \prod_{k=1}^{n} (1 - r_{kk}^{(k)}).$$

REMARK. As noted in [3] and [4], a direct consequence of the two theorems above is $\det(\text{Id} - u) \stackrel{\text{law}}{=} \prod_{k=1}^{n}(1 - r_{kk}^{(k)})$, where $u \sim \mu_{U(n)}$ and $r^{(k)} \sim \nu^{(k)}$, all being independent. As $r_{kk}^{(k)}$ is the first coordinate of a vector uniformly distributed on the unit complex sphere of dimension $n - k + 1$, it is equal in law to $e^{i\theta}\sqrt{B_{1,n-k}}$, with $\theta$ uniform on $(-\pi, \pi)$ and $B$ a beta variable with the indicated parameters. Consequently, with obvious notation for the following independent random variables,

$$(2.3) \qquad \det(\text{Id} - u) \stackrel{\text{law}}{=} \prod_{k=1}^{n}(1 - e^{i\theta_k}\sqrt{B_{1,n-k}}).$$

This relation will be useful in the proof of the Weyl integration formula, in the next section.

2.3. *Conditional Haar measure as the law of a product of independent reflections.* What could be the conditional expectation of $u \sim \mu_{U(n)}$, conditioned to have one eigenvalue at 1? As this conditioning is with respect to an event of measure 0, such a choice of conditional expectation is not trivial.

As previously, suppose we generate the Haar measure as a product of independent reflections: $u = r^{(1)} \cdots r^{(n)}$. Since $\mathrm{Id} - r^{(k)}$ has rank 1 a.s., our conditional expectation will naturally be constructed as a product of $n-1$ of these reflections: the unitary matrix $u$ has one eigenvalue $e^{i\theta} = 1$ if and only if $r^{(k)} = \mathrm{Id}$ for some $1 \leq k \leq n$, which yields $r_{kk}^{(k)} = 1$. As $r_{kk}^{(k)} \stackrel{\text{law}}{=} e^{i\theta}\sqrt{B_{1,n-k}}$, with the previous notation, $r_{nn}^{(n)}$ is more likely to be equal to 1 than any other $r_{kk}^{(k)}$ ($1 \leq k \leq n-1$).

Consequently, a good definition for the conditional expectation of $u \sim \mu_{U(n)}$, conditioned to have one eigenvalue at 1, is $r^{(1)} \cdots r^{(n-1)}$. This idea is formalized in the following way.

PROPOSITION 2.4. *Let $Z(x) = \det(x\,\mathrm{Id} - u)$ and $dx$ be the measure of $|Z(1)|$ under Haar measure on $U(n)$. There exists a continuous family of probability measures $P^{(x)}$ ($0 \leq x \leq 2^n$) such that for any Borel subset $\Gamma$ of $U(n)$*

$$(2.4) \qquad \mu_{U(n)}(\Gamma) = \int_0^{2^n} P^{(x)}(\Gamma)\, dx.$$

*Moreover $P^{(0)} = \nu^{(1)} \times \cdots \times \nu^{(n-1)}$ necessarily.*

REMARK. The continuity of the probability measures is in the sense of weak topology: the map

$$(2.5) \qquad x \mapsto \int_{U(n)} f(\omega)\, dP^{(x)}(\omega)$$

is continuous for any bounded continuous function $f$ on $U(n)$.

PROOF. We give an explicit expression of this conditional expectation, thanks to Theorem 2.3. Take $x > 0$. If $\prod_{k=1}^{n-1} |1 - r_{kk}^{(k)}| > x/2$, then there are two $r_{nn}^{(n)}$'s on the unit circle such that $\prod_{k=1}^{n} |1 - r_{kk}^{(k)}| = x$:

$$(2.6) \qquad r_{nn}^{(n)} = \exp\left(\pm 2i \arcsin \frac{x}{2 \prod_{k=1}^{n-1} |1 - r_{kk}^{(k)}|}\right).$$

These two numbers will be denoted $r_+$ et $r_-$. We write $\nu_\pm$ for the distribution of $r_\pm$, the random matrix in $\mathcal{R}^{(n)}$ equal to $\mathrm{Id}_{n-1} \oplus r_+$ with probability $1/2$, $\mathrm{Id}_{n-1} \oplus r_-$ with probability $1/2$. We define the conditional expectation, for any bounded continuous function $f$, by

$$(2.7) \quad \mathbb{E}_{\mu_{U(n)}}(f(u) \mid |Z(1)| = x) = \frac{\mathbb{E}(f(r^{(1)} \cdots r^{(n-1)} r_\pm) \mathbb{1}_{\prod_{k=1}^{n-1}|1-r_{kk}^{(k)}|>x/2})}{\mathbb{E}(\mathbb{1}_{\prod_{k=1}^{n-1}|1-r_{kk}^{(k)}|>x/2})},$$

the expectations on the RHS being with respect to $\nu^{(1)} \times \cdots \times \nu^{(n-1)} \times \nu_{\pm}$. For such a choice of the measures $P^{(x)}$ ($x > 0$), (2.4) holds thanks to Theorems 2.2 and 2.3. Moreover, these measures are continuous in $x$, and from (2.6) and (2.7) they converge to $\nu^{(1)} \times \cdots \times \nu^{(n-1)}$ as $x \to 0$. The continuity condition and formula (2.4) impose unicity for $(P^{(x)}, 0 \leq x \leq 2^n)$. Consequently $P^{(0)}$ necessarily coincides with $\nu^{(1)} \times \cdots \times \nu^{(n-1)}$. $\square$

For any $1 \leq k \leq n-1$, we can state some analogue of Proposition 2.4, conditioning now with respect to

$$(2.8) \qquad (|Z(1)|, |Z'(1)|, \ldots, |Z^{(k-1)}(1)|).$$

This leads to the following definition of the conditional expectation, which is the unique suitable choice preserving the continuity of measures with respect to (2.8).

DEFINITION 2.5. For any $1 \leq p \leq n-1$, $\nu^{(1)} \times \cdots \times \nu^{(n-p)}$ is called the Haar measure on $U(n)$ conditioned to have $p$ eigenvalues equal to 1.

REMARK. The above discussion can be held for the orthogonal group: the Haar measure on $O(n)$ conditioned to the existence of a stable subspace of dimension $p$ ($0 \leq p \leq n-1$) is

$$\nu_{\mathbb{R}}^{(1)} \times \cdots \times \nu_{\mathbb{R}}^{(n-p)},$$

where $\nu_{\mathbb{R}}^{(k)}$ is defined as the real analogue of $\nu^{(k)}$: a reflection $r$ is $\nu_{\mathbb{R}}^{(k)}$-distributed if $r(e_k)$ has its first $k-1$ coordinates equal to 0 and the others are uniformly distributed on the real unit sphere.

More generally, we can define this conditional Haar measure for any compact group generated by reflections, more precisely any compact group checking condition (R) in the sense of [4].

Take $p = n-1$ in Definition 2.5: the distribution of the unique eigenangle distinct from 1 coincides with the distribution of the nontrivial eigenangle of a reflection $r \sim \nu^{(1)}$, that is to say from (2.2)

$$e^{i\phi} = -\frac{1-r_{11}}{1-\overline{r_{11}}} \stackrel{\text{law}}{=} -\frac{1-e^{i\theta}\sqrt{B_{1,n-1}}}{1-e^{-i\theta}\sqrt{B_{1,n-1}}}.$$

In particular, this eigenvalue is not uniformly distributed on the unit circle: it converges in law to $-1$ as $n \to \infty$. This agrees with the idea of repulsion of the eigenvalues: we make it more explicit with the following probabilistic proof of the Weyl integration formula.

**3. A probabilistic proof of the Weyl integration formula.** The following two lemmas play a key role in our proof of the Weyl integration formula: the first shows that the spectral measure on $U(n)$ can be generated by $n-1$ reflections (instead of $n$) and the second one gives a transformation from this product of $n-1$ reflections in $U(n)$ to a product of $n-1$ reflections in $U(n-1)$, preserving the spectrum.

In the following, $u \stackrel{\text{sp}}{=} v$ means that the spectra of the matrices $u$ and $v$ are equally distributed.

3.1. *The conditioning lemma.* Remember that the measures $\nu^{(k)}$ ($1 \leq k \leq n$) are supported on the set of reflections: the following lemma would not be true by substituting our reflections with Householder transformations, for example.

LEMMA 3.1. *Take $r^{(k)} \sim \nu^{(k)}$ ($1 \leq k \leq n$), $\theta$ uniform on $(-\pi, \pi)$ and $u \sim \mu_{U(n)}$, all being independent. Then*

$$u \stackrel{\text{sp}}{=} e^{i\theta} r^{(1)} \cdots r^{(n-1)}.$$

PROOF. From Proposition 2.4, the spectrum of $r^{(1)} \cdots r^{(n-1)}$ is equal in law to the spectrum of $u$ conditioned to have one eigenvalue equal to 1.

Moreover, the Haar measure on $U(n)$ is invariant by translation, in particular by multiplication by $e^{i\phi}\operatorname{Id}$, for any fixed $\phi$: the distribution of the spectrum in invariant by rotation.

Consequently, the spectral distribution of $u \sim \mu_{U(n)}$ can be realized by successively conditioning to have one eigenvalue at 1 and then shifting by an independent uniform eigenangle, that is to say $u \stackrel{\text{sp}}{=} e^{i\theta} r^{(1)} \cdots r^{(n-1)}$, giving the desired result. □

3.2. *The slipping lemma.* Take $1 \leq k \leq n$ and $\delta$ a complex number. We first define a modification $\nu_\delta^{(k)}$ of the measure $\nu^{(k)}$ on the set of reflections $\mathcal{R}^{(k)}$. Let

$$\exp_\delta^{(k)} : \begin{cases} \mathcal{R}^{(k)} \to \mathbb{R}^+, \\ r \mapsto (1-r_{kk})^{\overline{\delta}}(1-\overline{r_{kk}})^{\delta}. \end{cases}$$

Then $\nu_\delta^{(k)}$ is defined as the $\exp_\delta^{(k)}$-sampling of a measure $\nu^{(k)}$ on $\mathcal{R}^{(k)}$, in the sense of the following definition.

DEFINITION 3.2. Let $(X, \mathcal{F}, \mu)$ be a probability space, and $h : X \mapsto \mathbb{R}^+$ a measurable function with $\mathbb{E}_\mu(h(x)) > 0$. Then a measure $\mu'$ is said to be the *h-sampling* of $\mu$ if for all bounded measurable functions $f$

$$\mathbb{E}_{\mu'}(f(x)) = \frac{\mathbb{E}_\mu(f(x)h(x))}{\mathbb{E}_\mu(h(x))}.$$

For $\mathfrak{Re}(\delta) > -1/2$, $0 < \mathbb{E}_{\nu^{(k)}}(\exp_\delta^{(k)}(r)) < \infty$, so $\nu_\delta^{(k)}$ is properly defined.

LEMMA 3.3. *Let $r^{(k)} \sim \nu^{(k)}$ ($1 \leq k \leq n-1$) and $r_1^{(k)} \sim \nu_1^{(k)}$ ($2 \leq k \leq n$) be $n \times n$ independent reflections. Then*

$$r^{(1)} \cdots r^{(n-1)} \stackrel{\mathrm{sp}}{=} r_1^{(2)} \cdots r_1^{(n)}.$$

PROOF. We proceed by induction on $n$. For $n = 2$, take $r \sim \nu^{(1)}$. Consider the unitary change of variables

$$(3.1) \quad \Phi : \begin{pmatrix} e_1 \\ e_2 \end{pmatrix} \mapsto \frac{1}{|1 - r_{11}|^2 + |r_{12}|^2} \begin{pmatrix} \overline{r_{12}} & -(1 - \overline{r_{11}}) \\ 1 - r_{11} & r_{12} \end{pmatrix} \begin{pmatrix} e_1 \\ e_2 \end{pmatrix}.$$

In this new basis, $r$ is diagonal with eigenvalues 1 and $r_{11} - |r_{12}|^2/(1 - \overline{r_{11}})$, so we only need to check that this last random variable is equal in law to the $|1 - X|^2$-sampling of a random variable $X$ uniform on the unit circle. This is a particular case of the identity in law given in Theorem A.1.

We now reproduce the above argument for general $n > 2$. Suppose the result is true at rank $n - 1$. Take $u \sim \mu_{U(n)}$, independent of all the other random variables. Obviously,

$$r^{(1)} \cdots r^{(n-1)} \stackrel{\mathrm{sp}}{=} (u^{-1} r^{(1)} u)(u^{-1} r^{(2)} \cdots r^{(n-1)} u).$$

As the uniform measure on the sphere is invariant by a unitary change of basis, by conditioning by $(u, r^{(2)}, \ldots, r^{(n-1)})$ we get

$$r^{(1)} \cdots r^{(n-1)} \stackrel{\mathrm{sp}}{=} r^{(1)}(u^{-1} r^{(2)} \cdots r^{(n-1)} u),$$

whose spectrum is equal in law (by induction) to the one of

$$r^{(1)}(u^{-1} r_1^{(3)} \cdots r_1^{(n)} u) \stackrel{\mathrm{sp}}{=} (u r^{(1)} u^{-1}) r_1^{(3)} \cdots r_1^{(n)} \stackrel{\mathrm{sp}}{=} r^{(1)} r_1^{(3)} \cdots r_1^{(n)}.$$

Consider now the change of basis $\Phi$ (3.1), extended to keep $(e_3, \ldots, e_n)$ invariant. As this transition matrix commutes with $r_1^{(3)} \cdots r_1^{(n)}$, to conclude we only need to show that $\Phi(r^{(1)}) \stackrel{\mathrm{law}}{=} r_1^{(2)}$. Both transformations are reflections, so a sufficient condition is $\Phi(r^{(1)})(e_2) \stackrel{\mathrm{law}}{=} r_1^{(2)}(e_2)$. A simple calculation gives

$$\Phi(r^{(1)})(e_2) = \left(0, r_{11} - \frac{|r_{12}|^2}{1 - \overline{r_{11}}}, c r_{13}, \ldots, c r_{1n}\right)^\top,$$

where the constant $c$ depends uniquely on $r_{11}$ and $r_{12}$. Hence the desired result is a direct consequence of the identity in law from Theorem A.1. □

REMARK. The above method and the identity in law stated in Theorem A.1 can be used to prove the following more general version of the slipping lemma. Let $1 \leq m \leq n-1$, $\delta_1, \ldots, \delta_m$ be complex numbers with real part greater than $-1/2$. Let $r_{\delta_k}^{(k)} \sim \nu_{\delta_k}^{(k)}$ ($1 \leq k \leq m$) and $r_{\delta_{k-1}+1}^{(k)} \sim \nu_{\delta_{k-1}+1}^{(k)}$ ($2 \leq k \leq m+1$) be independent $n \times n$ reflections. Then

$$r_{\delta_1}^{(1)} \cdots r_{\delta_m}^{(m)} \stackrel{\text{sp}}{=} r_{\delta_1+1}^{(2)} \cdots r_{\delta_m+1}^{(m+1)}.$$

In particular, iterating the above result,

(3.2) $$r^{(1)} \cdots r^{(n-p)} \stackrel{\text{sp}}{=} r_p^{(p+1)} \cdots r_p^{(n)}.$$

Together with Theorem 2.3, this implies that the eigenvalues of $r^{(1)} \cdots r^{(n-p)}$ have density (1.2), as expected.

3.3. *The proof by induction.* The two previous lemmas give a recursive proof of the following well-known result.

THEOREM 3.4. *Let $f$ be a class function on $U(n)$: $f(u) = f(\theta_1, \ldots, \theta_n)$, where the $\theta$'s are the eigenangles of $u$ and $f$ is symmetric. Then*

$$\mathbb{E}_{\mu_{U(n)}}(f(u)) = \frac{1}{n!} \int_{(-\pi,\pi)^n} f(\theta_1, \ldots, \theta_n) \prod_{1 \leq k < l \leq n} |e^{i\theta_k} - e^{i\theta_l}|^2 \frac{d\theta_1}{2\pi} \cdots \frac{d\theta_n}{2\pi}.$$

PROOF. We proceed by induction on $n$. The case $n=1$ is obvious. Suppose the result holds at rank $n-1$. Successively by the conditioning lemma and the slipping lemma, if $u \sim \mu_{U(n)}$,

$$u \stackrel{\text{sp}}{=} e^{i\theta} r^{(1)} \cdots r^{(n-1)} \stackrel{\text{sp}}{=} e^{i\theta} r_1^{(2)} \cdots r_1^{(n)}.$$

Hence, using the recurrence hypothesis, for any class function $f$,

$\mathbb{E}_{\mu_{U(n)}}(f(u))$

$$= \frac{1}{\text{cst}} \int_{-\pi}^{\pi} \frac{d\theta}{2\pi} \frac{1}{(n-1)!}$$

$$\times \int_{(-\pi,\pi)^{n-1}} f(\theta, \theta_2 + \theta, \ldots, \theta_n + \theta)$$

$$\times \prod_{2 \leq k < l \leq n} |e^{i\theta_k} - e^{i\theta_l}|^2 \prod_{j=2}^{n} |1 - e^{i\theta_j}|^2 \frac{d\theta_2}{2\pi} \cdots \frac{d\theta_n}{2\pi}$$

$$= \frac{1}{\text{cst}} \frac{1}{(n-1)!} \int_{(-\pi,\pi)^n} f(\theta_1, \ldots, \theta_n) \prod_{1 \leq k < l \leq n} |e^{i\theta_k} - e^{i\theta_l}|^2 \frac{d\theta_1}{2\pi} \cdots \frac{d\theta_n}{2\pi}.$$

Here cst comes from the sampling: $\text{cst} = \mathbb{E}_{\mu_{U(n-1)}}(|\det(\text{Id} - u)|^2)$. This is equal to $n$ from the decomposition of $\det(\text{Id} - u)$ into a product of independent random variables (2.3), which completes the proof. $\square$

**4. Derivatives as products of independent random variables.** As explained in [3] and [4], $\det(\mathrm{Id} - u)$ is equal in law to a product of $n$ independent random variables, for the Haar measure on $U(n)$ or $USp(2n)$, and $2n$ independent random variables, for the Haar measure on $SO(2n)$. These results are generalized to the Haar measures conditioned to the existence of a stable subspace with given dimension.

We first focus on the unitary group. Consider the conditional Haar measure (1.2) on $U(n)$ ($\theta_{n-p+1} = \cdots = \theta_n = 0$ a.s.). Then the $p$th derivative of the characteristic polynomial at 1 is

$$Z_U^{(p)} = p! \prod_{k=1}^{n-p} (1 - e^{i\theta_k}).$$

Our discussion about Haar measure on $U(n)$ conditioned to have a stable subspace with dimension $p$ allows us to decompose $Z_U^{(p)}$ as a product of independent random variables. More precisely, the following result holds.

COROLLARY 4.1. *Under the conditional Haar measure (1.2),*

$$\frac{Z_U^{(p)}}{p!} \stackrel{\mathrm{law}}{=} \prod_{\ell=1}^{n-p} (1 - X_\ell),$$

*where the $X_\ell$'s are independent random variables. The distribution of $X_\ell$ is the $|1-X|^{2p}$-sampling of a random variable $X = e^{i\theta}\sqrt{B_{1,\ell-1}}$, were $\theta$ is uniform on $(-\pi, \pi)$ and independently $B_{1,\ell-1}$ is a beta variable with the indicated parameters.*

PROOF. This proof directly relies on the suitable Definition 2.5 of the conditional Haar measure: it does not make use of the Weyl integration formula.

With the notation of Definition 2.5 ($r^{(1)} \cdots r^{(n-p)} \sim \nu^{(1)} \times \cdots \times \nu^{(n-p)}$),

$$Z_U^{(p)} \stackrel{\mathrm{law}}{=} \frac{d^p}{dx^p}\bigg|_{x=1} \det(x\,\mathrm{Id}_n - r^{(1)} \cdots r^{(n-p)}).$$

From (3.2), $r^{(1)} \cdots r^{(n-p)} \stackrel{\mathrm{sp}}{=} r_p^{(p+1)} \cdots r_p^{(n)}$, hence

$$Z_U^{(p)} \stackrel{\mathrm{law}}{=} \frac{d^p}{dx^p}\bigg|_{x=1} \det(x\,\mathrm{Id}_n - r_p^{(p+1)} \cdots r_p^{(n)}) = p! \prod_{k=p+1}^{n} (1 - \langle e_k, r_p^{(k)}(e_k)\rangle),$$

the last equality being a consequence of Theorem 2.3. The $r_p^{(k)}$'s are independent and $r_p^{(k)} \sim \nu_p^{(k)}$, which gives the desired result. $\square$

REMARK. Another proof of the above corollary consists in using the Weyl integration formula and the generalized Ewens sampling formula given in [4]. Actually, an identity similar to Corollary 4.1 can be stated for any Hua–Pickrell measure. More on these measures can be found in [2] for its connections with the theory of representations and in [4] for its analogies with the Ewens measures on permutation groups.

Corollary 4.1 admits an analogue for the Jacobi ensemble on the segment. Indeed, Lemma 5.2 and Proposition 5.3 in [9], by Killip and Nenciu, immediately imply that under the probability measure (cst is the normalization constant)

$$\text{(4.1)} \qquad \text{cst}|\Delta(x_1,\ldots,x_n)|^\beta \prod_{j=1}^{n}(2-x_j)^a(2+x_j)^b\, dx_1\cdots dx_n$$

on $(-2,2)^n$, the following identity in law holds (the $x_k$'s being the eigenvalues of a matrix $u$):

$$\text{(4.2)} \quad (\det(2\operatorname{Id}-u), \det(2\operatorname{Id}+u)) \stackrel{\text{law}}{=} \left(2\prod_{k=0}^{2n-2}(1-\alpha_k), 2\prod_{k=0}^{2n-2}(1+(-1)^k\alpha_k)\right),$$

with the $\alpha_k$'s independent with density $f_{s(k),t(k)}$ ($f_{s,t}$ is defined below) on $(-1,1)$ with

$$\begin{cases} s(k) = \dfrac{2n-k-2}{4}\beta + a + 1, \\ t(k) = \dfrac{2n-k-2}{4}\beta + b + 1, & \text{if } k \text{ is even,} \\ s(k) = \dfrac{2n-k-3}{4}\beta + a + b + 2, \\ t(k) = \dfrac{2n-k-1}{4}\beta, & \text{if } k \text{ is odd.} \end{cases}$$

DEFINITION 4.2. The density $f_{s,t}$ on $(-1,1)$ is

$$f_{s,t}(x) = \frac{2^{1-s-t}\Gamma(s+t)}{\Gamma(s)\Gamma(t)}(1-x)^{s-1}(1+x)^{t-1}.$$

Moreover, for $X$ with density $f_{s,t}$, $\mathbb{E}(X) = \frac{t-s}{t+s}$, $\mathbb{E}(X^2) = \frac{(t-s)^2+(t+s)}{(t+s)(t+s+1)}$.

Consequently, the analogue of Corollary 4.1 can be stated for $SO(2n)$ and $USp(2n)$, relying on formula (4.2). More precisely, the $2n$ eigenvalues of $u \in SO(2n)$ or $USp(2n)$ are pairwise conjugated, and noted $(e^{\pm i\theta_1},\ldots,e^{\pm i\theta_n})$. The Weyl integration formula states that the eigenvalues statistics are

$$\text{cst} \prod_{1\leq k<\ell\leq n}(\cos\theta_k - \cos\theta_\ell)^2\, d\theta_1\cdots d\theta_n$$

on $SO(2n)$. On the symplectic group $USp(2n)$, these statistics are

$$\text{cst} \prod_{1\leq k<\ell\leq n}(\cos\theta_k - \cos\theta_\ell)^2 \prod_{i=1}^{n}(1-\cos\theta_i)(1+\cos\theta_i)\,d\theta_1\cdots d\theta_n.$$

Hence, the change of variables

$$x_i = 2\cos\theta_j$$

implies the following link between $SO(2n)$, $USp(2n)$ and the Jacobi ensemble:

- On $SO(2n+2p^+ +2p^-)$, endowed with its Haar measure conditioned to have $2p^+$ eigenvalues at 1 and $2p^-$ at $-1$, the distribution of $(x_1,\ldots,x_n)$ is the Jacobi ensemble (4.1) with parameters $\beta = 2$, $a = 2p^+ - \frac{1}{2}$, $b = 2p^- - \frac{1}{2}$.
- On $USp(2n+2p^{(+)}+2p^-)$, endowed with its Haar measure conditioned to have $2p^+$ eigenvalues at 1 and $2p^-$ at $-1$, the distribution of $(x_1,\ldots,x_n)$ is the Jacobi ensemble (4.1) with parameters $\beta = 2$, $a = 2p^+ + \frac{1}{2}$, $b = 2p^- + \frac{1}{2}$.

Moreover, for the above groups $\mathcal{G} = SO(2n+2p^+ +2p^-)$ or $Sp(2n+2p^+ + 2p^-)$ with $2p^+$ eigenvalues at 1 and $2p^-$ at $-1$, $Z_\mathcal{G}^{(2p^+)}$ denotes the $2p^+$th derivative of the characteristic polynomial at point 1 and $Z_\mathcal{G}^{(2p^-)}$ the $2p^-$th derivative of the characteristic polynomial at point $-1$:

$$\begin{cases} \dfrac{Z_\mathcal{G}^{(2p^+)}}{(2p^+)!2^{p^-}} = \prod_{k=1}^{n}(1-e^{i\theta_k})(1-e^{-i\theta_k}) = \prod_{k=1}^{n}(2-x_k), \\ \dfrac{Z_\mathcal{G}^{(2p^-)}}{(2p^-)!2^{p^+}} = \prod_{k=1}^{n}(-1-e^{i\theta_k})(-1-e^{-i\theta_k}) = \prod_{k=1}^{n}(2+x_k). \end{cases}$$

Combining this with formula (4.2) leads to the following analogue of Corollary 4.1.

COROLLARY 4.3. *With the above notation and definition of conditional spectral Haar measures on $SO(2n+2p^+ +2p^-)$,*

$$\left(\frac{Z_{SO}^{(2p^+)}}{(2p^+)!2^{p^-}}, \frac{Z_{SO}^{(2p^-)}}{(2p^-)!2^{p^+}}\right) \stackrel{\text{law}}{=} \left(2\prod_{k=0}^{2n-2}(1-X_k), 2\prod_{k=0}^{2n-2}(1+(-1)^k X_k)\right),$$

where the $X_k$'s are independent and $X_k$ with density $f_{s(k),t(k)}$ on $(-1,1)$ given by Definition 4.2 with parameters

$$\begin{cases} s(k) = \dfrac{2n-k-1}{2} + 2p^+, \\ t(k) = \dfrac{2n-k-1}{2} + 2p^-, & \text{if } k \text{ is even,} \\ s(k) = \dfrac{2n-k-1}{2} + 2p^+ + 2p^-, \\ t(k) = \dfrac{2n-k-1}{2}, & \text{if } k \text{ is odd.} \end{cases}$$

The same result holds for the joint law of $Z_{USp}^{(2p^+)}$ and $Z_{USp}^{(2p^-)}$, but with the parameters

$$\begin{cases} s(k) = \dfrac{2n-k+1}{2} + 2p^+, \\ t(k) = \dfrac{2n-k+1}{2} + 2p^-, & \text{if } k \text{ is even,} \\ s(k) = \dfrac{2n-k+3}{2} + 2p^+ + 2p^-, \\ t(k) = \dfrac{2n-k-1}{2}, & \text{if } k \text{ is odd.} \end{cases}$$

**5. Limit theorems for the derivatives.** Let $Z_{SO}^{(2p)}$ be the $2p$th derivative of the characteristic polynomial at point 1, for the Haar measure on $SO(n+2p)$ conditioned to have $2p$ eigenvalues equal to 1. In the study of moments of $L$-functions associated to elliptic curves, Snaith explains that the moments of $Z_{SO}^{(2p)}$ are relevant: she conjectures that $Z_{SO}^{(2p)}$ is related to averages on $L$-functions moments and therefore, via the Swinnerton–Dyer conjecture, on the rank of elliptic curves. For the number theoretic applications of these derivatives, see [10, 13] and [14].

Relying on the Selberg integral, she computed the asymptotics of the density of $Z_{SO}^{(2p)}$ as $\varepsilon \to 0$, finding

$$(5.1) \qquad \mathbb{P}(Z_{SO}^{(2p)} < \varepsilon) \overset{\varepsilon \to 0}{\sim} c_{n,p} \varepsilon^{2p+1/2},$$

for an explicit constant $c_{n,p}$. Similar results (and also central limit theorems) are given in this section for the symplectic and unitary groups.

5.1. *Limit densities.* Let $(x_1, \ldots, x_n)$ have the Jacobi distribution (4.1) on $(-2,2)^n$. The asymptotics of the density of

$$\det^{(+)} := \prod_{k=1}^{n}(2-x_i)$$

near 0 can be easily evaluated from (4.2). Indeed, let $f$ be a continuous function and $h_n$ denote the density of $\det^{(+)}$ on $(0, \infty)$. With the notation of (4.2), as $\alpha_{2n-2}$ has law $f_{a+1,b+1}$,

$$\mathbb{E}(f(\det^{(+)})) = c \int_{-1}^{1} (1-x)^a (1+x)^b \mathbb{E}\left( f\left( 2(1-x) \prod_{k=0}^{2n-3} (1-\alpha_k) \right) \right) dx$$

with $c = 2^{-1-a-b} \Gamma(a+b+2)/(\Gamma(a+1)\Gamma(b+1))$. The change of variable $\varepsilon = 2(1-x) \prod_{k=0}^{2n-3}(1-\alpha_k)$ therefore yields

$$h_n(\varepsilon) = c \mathbb{E}\left( \left( \frac{1}{2\prod_{k=0}^{2n-3}(1-\alpha_k)} \right)^{a+1} \left( 2 - \frac{\varepsilon}{2\prod_{k=0}^{2n-3}(1-\alpha_k)} \right)^b \right) \varepsilon^a,$$

implying immediately the following corollary of Killip and Nenciu's formula (4.2).

COROLLARY 5.1. *For the Jacobi distribution (4.1) on $(-2,2)^n$, the density of the characteristic polynomial $\det^{(+)}$ near 0 is, for some constant $c(n)$,*

$$h_n(\varepsilon) \stackrel{\varepsilon \to 0}{\sim} c(n) \varepsilon^a.$$

Note that this constant is effective:

$$c(n) = \frac{\Gamma(a+b+2)}{2^{2(1+a)} \Gamma(a+1) \Gamma(b+1)} \prod_{k=0}^{2n-3} \mathbb{E}\left( \left( \frac{1}{1-\alpha_k} \right)^{1+a} \right).$$

As an application of Corollary 5.1, the correspondence $a = 2p - \frac{1}{2}$ shows that for the Haar measure on $SO(2n+2p)$, conditionally to have $2p$ eigenvalues equal to 1, this density has order $\varepsilon^{2p-1/2}$, which agrees with (5.1). Of course, Corollary 5.1 gives the same way the asymptotic density of the characteristic polynomial for the symplectic $(a = 2p + 1/2)$ groups or the orthogonal groups with odd dimension.

Moreover, the same method, based on Corollary 4.1, gives an analogous result for the unitary group.

COROLLARY 5.2. *Let $h_n^{(U)}$ be the density of $|Z_U^{(p)}|$, with the notation of Corollary 4.1. Then, for some constant $d(n)$,*

$$h_n^{(U)}(\varepsilon) \stackrel{\varepsilon \to 0}{\sim} d(n) \varepsilon^{2p}.$$

REMARK. With a similar method [decomposition of $\det^{(+)}$ as a product of independent random variables], such asymptotics were already obtained by Yor [15] for the density of the characteristic polynomial on the group $SO(n)$.

5.2. *Central limit theorems.* From (4.2), $\log\det(2\,\mathrm{Id}-u)$ and $\log\det(2\,\mathrm{Id}+u)$ [resp., abbreviated as $\log\det^{(+)}$ and $\log\det^{(-)}$] can be jointly decomposed as sums of independent random variables. Hence, the classical central limit theorems in probability theory imply the following result. Note that, despite the correlation appearing from (4.2), $\log\det^{(+)}$ and $\log\det^{(-)}$ are independent in the limit.

THEOREM 5.3. *Let $u$ have spectral measure the Jacobi ensemble (4.1), with $\beta > 0$, $a,b \geq 0$. Then*

$$\left(\frac{\log\det^{(+)} + (1/2 - (2a+1)/\beta)\log n}{\sqrt{2/\beta \log n}}, \frac{\log\det^{(-)} + (1/2 - (2b+1)/\beta)\log n}{\sqrt{2/\beta \log n}}\right)$$

$$\xrightarrow{\mathrm{law}} (\mathcal{N}_1, \mathcal{N}_2)$$

*as $n \to \infty$, with $\mathcal{N}_1$ and $\mathcal{N}_2$ independent standard normal variables.*

PROOF. We keep the notation from (4.2):

$$\begin{cases} \log\det^{(+)} = \log 2 + \sum_{\text{odd } k} \log(1-\alpha_k) + \sum_{\text{even } k} \log(1-\alpha_k), \\ \log\det^{(-)} = \log 2 + \sum_{\text{even } k} \log(1-\alpha_k) + \sum_{\text{even } k} \log(1+\alpha_k), \end{cases}$$

with $0 \leq k \leq 2n-2$. Let us first consider $X_n = \sum_{\text{odd } k} \log(1-\alpha_k)$. From (4.2), $X_n \stackrel{\mathrm{law}}{=} \sum_{k=1}^{n-1} \log(1-x_k)$ with independent $x_k$'s, $x_k$ having density $f_{(k-1)\beta/2+a+b+2,k\beta/2}$. In particular, $\mathbb{E}(x_k) = \frac{-a-b-2+\beta/2}{\beta k} + O(\frac{1}{k^2})$ and $\mathbb{E}(x_k^2) = \frac{1}{\beta k} + O(\frac{1}{k^2})$. From the Taylor expansion of $\log(1-x)$,

$$X_n \stackrel{\mathrm{law}}{=} \underbrace{\sum_{k=1}^{n-1} \left(-x_k - \frac{x_k^2}{2}\right)}_{X_n^{(1)}} - \underbrace{\sum_{k=1}^{n-1} \sum_{\ell \geq 3} \frac{x_k^\ell}{\ell}}_{X_n^{(2)}}.$$

Let $X = \sum_{k=1}^{\infty} \sum_{\ell \geq 3} \frac{|x_k|^\ell}{l}$. A calculation implies $\mathbb{E}(X) < \infty$, so $X < \infty$ a.s. and consequently, $|X_n^{(2)}|/\sqrt{\log n} \leq X/\sqrt{\log n} \to 0$ a.s. as $n \to \infty$. Moreover,

$$\mathbb{E}\left(-x_k - \frac{x_k^2}{2}\right) = \frac{a+b+3/2-\beta/2}{\beta k} + O\left(\frac{1}{k^2}\right),$$

$$\mathrm{Var}\left(-x_k - \frac{x_k^2}{2}\right) = \frac{1}{\beta k} + O\left(\frac{1}{k^2}\right),$$

so the classical central limit theorem (see, e.g., [11]) implies that $(X_n^{(1)} - \frac{a+b+3/2-\beta/2}{\beta}\log n)/\sqrt{(\log n)/\beta} \xrightarrow{\mathrm{law}} \mathcal{N}_1$ as $n \to \infty$, with $\mathcal{N}_1$ a standard normal random variable. Gathering the convergences for $X_n^{(1)}$ and $X_n^{(2)}$ gives

$$\text{(5.2)} \qquad \frac{X_n - (a+b+3/2-\beta/2)/\beta \log n}{\sqrt{1/\beta \log n}} \xrightarrow{\text{law}} \mathcal{N}_1.$$

We now concentrate on $Y_n^{(+)} = \sum_{\text{even } k} \log(1-\alpha_k)$ and $Y_n^{(-)} = \sum_{\text{even } k} \log(1+\alpha_k)$. From (4.2) $(Y_n^{(+)}, Y_n^{(-)}) \stackrel{\text{law}}{=} (\sum_1^n \log(1-y_k), \sum_1^n \log(1+y_k))$ with independent $y_k$'s, $y_k$ having density $f_{(k-1)\beta/2+a+1,(k-1)\beta/2+b+1}$. We now have

$$\mathbb{E}\left(\pm y_k - \frac{y_k^2}{2}\right) = \frac{\pm(b-a) - 1/2}{\beta k} + O\left(\frac{1}{k^2}\right),$$

$$\text{Var}\left(\pm y_k - \frac{y_k^2}{2}\right) = \frac{1}{\beta k} + O\left(\frac{1}{k^2}\right).$$

Consequently, as previously the two first terms in the Taylor expansions of $\log(1 \pm y_k)$ can be isolated to get the following central limit theorem for any real numbers $\lambda^{(+)}$ and $\lambda^{(-)}$:

$$\lambda^{(+)} \frac{Y_n^{(+)} - (a-b-1/2)/\beta \log n}{\sqrt{1/\beta \log n}} + \lambda^{(-)} \frac{Y_n^{(-)} - (b-a-1/2)/\beta \log n}{\sqrt{1/\beta \log n}}$$

$$\text{(5.3)} \qquad \xrightarrow{\text{law}} (\lambda^{(+)} - \lambda^{(-)})\mathcal{N}_2$$

with $\mathcal{N}_2$ a standard normal variable, independent of $\mathcal{N}_1$, because the odd and even $\alpha_k$'s are independent. Gathering convergences (5.2) and (5.3) shows that $\left(\frac{\log \det^{(+)} + (1/2-(2a+1)/\beta) \log n}{\sqrt{2/\beta \log n}}, \frac{\log \det^{(-)} + (1/2-(2b+1)/\beta) \log n}{\sqrt{2/\beta \log n}}\right)$ converges in law to $\frac{1}{\sqrt{2}}(\mathcal{N}_1 + \mathcal{N}_2, \mathcal{N}_1 - \mathcal{N}_2)$. The covariance matrix of this Gaussian vector is diagonal, hence its coordinates are independent, concluding the proof. $\square$

An immediate corollary of the previous theorem concerns the derivatives of characteristic polynomials on $SO(2n)$ and $USp(2n)$.

COROLLARY 5.4. *With the notation of Corollary 4.3,*

$$\left(\frac{\log Z_{SO}^{(2p^+)} - (2p^+ - 1/2) \log n}{\sqrt{\log n}}, \frac{\log Z_{SO}^{(2p^-)} - (2p^- - 1/2) \log n}{\sqrt{\log n}}\right)$$

$$\xrightarrow{\text{law}} (\mathcal{N}_1, \mathcal{N}_2)$$

*as $n \to \infty$, with $\mathcal{N}_1$ and $\mathcal{N}_2$ independent standard normal variables. The same result holds on the symplectic group conditioned to have $2p^+$ eigenvalues at 1 and $2p^-$ at $-1$, but with the parameters $2p^{(+)}$ and $2p^{(-)}$ replaced by $2p^{(+)}+1$ and $2p^{(-)}+1$ in the above formula.*

These central limit theorems about orthogonal and symplectic groups have an analogue for the unitary group. We only state it, the proof being very similar to the previous one and relying on the decomposition as a product of independent random variables, Corollary 4.1. In the following the complex logarithm is defined continuously on $(0,1)$ as the value of $\log\det(\mathrm{Id} - xu)$ from $x = 0$ to $x = 1$.

COROLLARY 5.5. *If $Z_U^{(p)}$ is the pth derivative of the characteristic polynomial at 1, for the Haar measure on $U(n)$ conditioned to have $p$ eigenvalues equal to 1,*

$$\frac{\log Z_U^{(p)} - p \log n}{\sqrt{\log n}} \xrightarrow{\text{law}} \frac{1}{\sqrt{2}}(\mathcal{N}_1 + i\mathcal{N}_2)$$

*as $n \to \infty$, $\mathcal{N}_1$ and $\mathcal{N}_2$ being independent standard normal variables.*

## APPENDIX: AN EQUALITY IN LAW

**Known results about beta variables.** A beta random variable $B_{a,b}$ with strictly positive coefficients $a$ and $b$ has distribution on $(0,1)$ given by

$$\mathbb{P}\{B_{a,b} \in dt\} = \frac{\Gamma(a+b)}{\Gamma(a)\Gamma(b)} t^{a-1}(1-t)^{b-1}\,dt.$$

Its Mellin transform is

(A.1) $$\mathbb{E}[B_{a,b}^s] = \frac{\Gamma(a+s)}{\Gamma(a)} \frac{\Gamma(a+b)}{\Gamma(a+b+s)}, \qquad s > 0.$$

For the unform measure on the real sphere

$$\mathscr{S}_\mathbb{R}^n = \{(r_1, \ldots, r_n) \in \mathbb{R}^n : r_1^2 + \cdots + r_n^2 = 1\},$$

the sum of the squares of the first $k$ coordinates is equal in law to $B_{k/2,(n-k)/2}$. Consequently, under the uniform measure on the complex sphere $\mathscr{S}_\mathbb{C}^n = \{(c_1, \ldots, c_n) \in \mathbb{C}^n : |c_1|^2 + \cdots + |c_n|^2 = 1\}$,

$$c_1 \stackrel{\text{law}}{=} e^{i\theta}\sqrt{B_{1,n-1}}$$

with $\theta$ uniform on $(-\pi, \pi)$ and independent from $B_{1,n-1}$.

**The identity.** Using Mellin–Fourier transforms, we will prove the following equality in law.

THEOREM A.1. *Let $\mathfrak{Re}(\delta) > -1/2$, $\lambda > 1$. Take independently:*

- *$\theta$ uniform on $(-\pi, \pi)$ and $B_{1,\lambda}$ a beta variable with the indicated parameters; $Y$ distributed as the $(1-x)^{\overline{\delta}}(1-\overline{x})^\delta$-sampling of $x = e^{i\theta_1}\sqrt{B_{1,\lambda}}$;*

- $B_{1,\lambda-1}$ be a beta variable with the indicated parameters;
- $Z$ distributed as the $(1-x)^{\overline{\delta}+1}(1-\overline{x})^{\delta+1}$-sampling of $x = e^{i\theta}\sqrt{B_{1,\lambda-1}}$.

Then

$$Y - \frac{(1-|Y|^2)B_{1,\lambda-1}}{1-\overline{Y}} \stackrel{\text{law}}{=} Z.$$

PROOF. Actually, we will show that

$$X = 1 - \left(Y - \frac{(1-|Y|^2)B_{1,\lambda-1}}{1-\overline{Y}}\right) \stackrel{\text{law}}{=} 1 - Z.$$

First note that, by Lemma A.4,

$$(A.2) \qquad 1 - Y \stackrel{\text{law}}{=} 2\cos\phi\, e^{i\phi} B_{\lambda+\delta+\overline{\delta}+1,\lambda},$$

where $\phi$ has probability density $c(1+e^{2i\phi})^{\lambda+\delta}(1+e^{-2i\phi})^{\lambda+\overline{\delta}}\mathbb{1}_{(-\pi/2,\pi/2)}$ ($c$ is the normalization constant). Consequently, by a straightforward calculation,

$$X \stackrel{\text{law}}{=} 2\cos\phi\, e^{i\phi}(B_{\lambda+\delta+\overline{\delta}+1,\lambda} + (1-B_{\lambda+\delta+\overline{\delta}+1,\lambda})B_{1,\lambda-1}).$$

Consider the uniform distribution on $\mathscr{S}_{\mathbb{R}}^{2(2\lambda+\delta+\overline{\delta}+1)-1}$. Then the sum of the squares of the first $2(\lambda+\delta+\overline{\delta}+2)$ coordinates is equal in law to $B_{\lambda+\delta+\overline{\delta}+2,\lambda-1}$, but also to $B_{\lambda+\delta+\overline{\delta}+1,\lambda}+(1-B_{\lambda+\delta+\overline{\delta}+1,\lambda})B_{1,\lambda-1}$ by counting the first $2(\lambda+\delta+\overline{\delta}+1)$ coordinates first and then the next two. Hence

$$X \stackrel{\text{law}}{=} 2\cos\phi\, e^{i\phi} B_{\lambda+\delta+\overline{\delta}+2,\lambda-1}.$$

Consequently Lemma A.3 implies the following Mellin–Fourier transform

$$\mathbb{E}(|X|^t e^{is\arg X}) = \frac{\Gamma(\lambda+\delta+1)\Gamma(\lambda+\overline{\delta}+1)}{\Gamma(\lambda+\delta+\overline{\delta}+2)}$$
$$\times \frac{\Gamma(\lambda+\delta+\overline{\delta}+2+t)}{\Gamma(\lambda+\delta+(t+s)/2+1)\Gamma(\lambda+\overline{\delta}+(t-s)/2+1)}.$$

Using Lemma A.2, the Mellin–Fourier transform of $1-Z$, coincides with the above expression, completing the proof. $\square$

LEMMA A.2. *Let $\lambda > 0$ and $X = 1 + e^{i\theta}\sqrt{B_{1,\lambda}}$, where $\theta$, uniformly distributed on $(-\pi,\pi)$, is assumed independent from $B_{1,\lambda}$. Then, for all $t$ and $s$ with $\mathfrak{Re}(t\pm s) > -1$*

$$\mathbb{E}(|X|^t e^{is\arg X}) = \frac{\Gamma(\lambda+1)\Gamma(\lambda+1+t)}{\Gamma(\lambda+1+(t+s)/2)\Gamma(\lambda+1+(t-s)/2)}.$$

PROOF. First, note that

$$\mathbb{E}(|X|^t e^{is \arg X}) = \mathbb{E}((1 + e^{i\theta}\sqrt{B_{1,\lambda}})^a (1 + e^{-i\theta}\sqrt{B_{1,\lambda}})^b)$$

with $a = (t+s)/2$ and $b = (t-s)/2$. Recall that if $|x| < 1$ and $u \in \mathbb{R}$ then

$$(1+x)^u = \sum_{k=0}^{\infty} \frac{(-1)^k (-u)_k}{k!} x^k,$$

where $(y)_k = y(y+1)\cdots(y+k-1)$ is the Pochhammer symbol. As $|e^{i\theta}\sqrt{B_{1,\lambda}}| < 1$ a.s., we get

$$\mathbb{E}[|X|^t e^{is \arg X}] = \mathbb{E}\left[\left(\sum_{k=0}^{\infty} \frac{(-1)^k (-a)_k}{k!} B_{1,\lambda}^{k/2} e^{ik\theta}\right)\left(\sum_{\ell=0}^{\infty} \frac{(-1)^\ell (-b)_\ell}{\ell!} B_{1,\lambda}^{\ell/2} e^{-i\ell\theta}\right)\right].$$

After expanding this double sum [it is absolutely convergent because $\mathfrak{Re}(t \pm s) > -1$], all terms with $k \neq \ell$ will give an expectation equal to 0, so as $\mathbb{E}(B_{1,\lambda}^k) = \frac{\Gamma(1+k)\Gamma(\lambda+1)}{\Gamma(1)\Gamma(\lambda+1+k)} = \frac{k!}{(\lambda+1)_k}$,

$$\mathbb{E}[|X|^t e^{is \arg X}] = \sum_{k=0}^{\infty} \frac{(-a)_k(-b)_k}{k!(\lambda+1)_k}.$$

Note that this series is equal to the value at $z=1$ of the hypergeometric function $_2F_1(-a, -b, \lambda+1; z)$. This value is well known (see, e.g., [1]) and yields:

$$\mathbb{E}[|X|^t e^{is \arg X}] = \frac{\Gamma(\lambda+1)\Gamma(\lambda+1+a+b)}{\Gamma(\lambda+1+a)\Gamma(\lambda+1+b)}.$$

This is the desired result. □

LEMMA A.3. *Take $\phi$ with probability density $c(1+e^{2i\phi})^{\overline{z}}(1+e^{-2i\phi})^z \times \mathbb{1}_{(-\pi/2, \pi/2)}$, where $c$ is the normalization constant, $\mathfrak{Re}(z) > -1/2$. Let $X = 2\cos\phi \, e^{i\phi}$. Then*

(A.3)
$$\mathbb{E}[|X|^t e^{is \arg X}] = \frac{\Gamma(z+1)\Gamma(\overline{z}+1)}{\Gamma(z+\overline{z}+1)}$$
$$\times \frac{\Gamma(z+\overline{z}+t+1)}{\Gamma(\overline{z}+(t+s)/2+1)\Gamma(z+(t-s)/2+1)}.$$

PROOF. From the definition of $X$,

$$\mathbb{E}(|X|^t e^{is \arg X}) = c \int_{-\pi/2}^{\pi/2} (1+e^{2ix})^{\overline{z}+(t+s)/2}(1+e^{-2ix})^{z+(t-s)/2} \, dx.$$

Both terms on the RHS can be expanded as a series in $e^{2ix}$ or $e^{-2ix}$ for all $x \neq 0$. Integrating over $x$ between $-\pi/2$ and $\pi/2$, only the diagonal terms remain, and so

$$\mathbb{E}(|X|^t e^{is \arg X}) = c \sum_{k=0}^{\infty} \frac{(-\overline{z} - (+t+s)/2)_k (-z - (t-s)/2)_k}{(k!)^2}.$$

The value of a $_2F_1$ hypergeometric functions at $z = 1$ is well known (see, e.g., [1]), hence

$$\mathbb{E}(|X|^t e^{is \arg X}) = c \frac{\Gamma(z + \overline{z} + t + 1)}{\Gamma(\overline{z} + (t+s)/2 + 1)\Gamma(z + (t-s)/2 + 1)}.$$

As this is 1 when $s = t = 0$, $c = \frac{\Gamma(z+1)\Gamma(\overline{z}+1)}{\Gamma(z+\overline{z}+1)}$, which completes the proof. $\square$

LEMMA A.4. *Let $\lambda > 2$, $\mathfrak{Re}(\delta) > -1/2$, $\theta$ uniform on $(-\pi, \pi)$ and $B_{1,\lambda-1}$ a beta variable with the indicated parameters. Let $Y$ be distributed as the $(1-\overline{x})^\delta (1-x)^{\overline{\delta}}$-sampling of $x = e^{i\theta} \sqrt{B_{1,\lambda-1}}$. Then*

$$1 - Y \stackrel{\text{law}}{=} 2 \cos \phi e^{i\phi} B_{\lambda+\delta+\overline{\delta}, \lambda-1}$$

*with $B_{\lambda+\delta+\overline{\delta}, \lambda-1}$ a beta variable with the indicated parameters and, independently, $\phi$ having probability density $c(1+e^{2i\phi})^{\lambda+\delta-1}(1+e^{-2i\phi})^{\lambda+\overline{\delta}-1}\mathbb{1}_{(-\pi/2,\pi/2)}$ ($c$ is the normalization constant).*

PROOF. The Mellin–Fourier transform of $X = 1 - Y$ can be evaluated using Lemma A.2, and equals

$$\mathbb{E}(|X|^t e^{is \arg X}) = \frac{\Gamma(\lambda + \delta)\Gamma(\lambda + \overline{\delta})\Gamma(\lambda + t + \delta + \overline{\delta})}{\Gamma(\lambda + (t+s)/2 + \delta)\Gamma(\lambda + (t-s)/2 + \overline{\delta})\Gamma(\lambda + \delta + \overline{\delta})}.$$

On the other hand, using Lemma A.3 and (A.1), the Mellin–Fourier transform of $2 \cos \phi e^{i\phi} B_{n+\delta+\overline{\delta}, \lambda}$ coincides with the above result. $\square$

TELECOM PARIS TECH
46 RUE BARRAULT
75634 PARIS CEDEX 13
FRANCE
E-MAIL: bourgade@enst.fr